\documentclass[12pt,a4paper,reqno]{amsart}
\usepackage{amsmath}
\usepackage{amsfonts}
\usepackage{amssymb}
\numberwithin{equation}{section}

\theoremstyle{plain}

\newtheorem{theorem}[subsection]{Theorem}
\newtheorem{proposition}[subsection]{Proposition}
\newtheorem{lemma}[subsection]{Lemma}

\theoremstyle{definition}

\newtheorem{definition}[subsection]{Definition}

\newtheorem{question}[subsection]{Question}

\renewcommand{\leq}{\leqslant}
\renewcommand{\geq}{\geqslant}

\newsavebox{\proofbox}
\savebox{\proofbox}{\begin{picture}(7,7)%
  \put(0,0){\framebox(7,7){}}\end{picture}}
\def\boxeq{\tag*{\usebox{\proofbox}}}

\newcommand{\md}[1]{\ensuremath{(\mbox{mod}\, #1)}}
\newcommand{\mdsub}[1]{\ensuremath{(\mbox{\scriptsize mod}\, #1)}}

\def\proof{\noindent\textit{Proof. }}
\def\endproof{\hfill{\usebox{\proofbox}}}
\def\stp{\mbox{\texttt{STEP}}}
\def\E{\mathbb{E}}
\def\Z{\mathbb{Z}}

\def\R{\mathbb{R}}

\def\Zp{\mathbb{Z}/p\mathbb{Z}}

\def\ni{\noindent}
\def\vs{\vspace{11pt}}

\begin{document}

\title{On the Littlewood problem modulo a prime}

\author{Ben Green}
\address{Department of Mathematics\\
University of Bristol\\
University Walk\\
Bristol BS8 1TW\\
England
}
\email{b.j.green@bristol.ac.uk}

\author{Sergei Konyagin}
\address{Department of Mathematics and Mechanics\\
Moscow State University\\
Moscow 119992, Russia
}
\email{konyagin@ok.ru}

\thanks{This research was started when the authors visited PIMS, Vancouver, in Spring 2004. It was partially conducted during the period the first-named author served as a Clay Research Fellow. He would like to express his sincere gratitude to the Clay In
stitute, and also to the Massachusetts Institute of Technology, where he was a
Visiting Professor for the academic year 2005-06. The second-named author was
supported by the INTAS grant 03-51-5070}

\begin{abstract}
\noindent Let $p$ be a prime, and let $f : \mathbb{Z}/p\mathbb{Z} \rightarrow \mathbb{R}$ be a function with $\E f = 0$ and $\Vert \widehat{f} \Vert_1 \leq 1$. Then \[ \min_{x \in \Zp} |f(x)| = O(\log p)^{-1/3 + \epsilon}.\]
One should think of $f$ as being ``approximately continuous''; our result is then an ``approximate intermediate value theorem''.\vs

\ni As an immediate consequence we show that if $A \subseteq \Zp$ is a set of cardinality $\lfloor p/2\rfloor$ then
$\sum_r |\widehat{1_A}(r)| \gg (\log p)^{1/3 - \epsilon}$. This gives a result on a ``mod $p$'' analogue of Littlewood's well-known problem concerning the smallest possible $L^1$-norm of the Fourier transform of a set of $n$ integers.\vs

\ni Another application is to answer a question of Gowers. If $A \subseteq \Zp$ is a set of size $\lfloor p/2 \rfloor$ then there is some $x \in \Zp$ such that
\[ \left| |A \cap (A + x)| - p/4 \right| = o(p).\]
 \end{abstract}

\maketitle

\section{Introduction} \label{sec1}

\noindent The problem of Littlewood to which the title refers is the following question. If $A \subseteq \mathbb{Z}$ is a set of $n$ integers, what is the smallest possible value of the $L^1$-norm of the exponential sum over $A$. That is to say, what is
\[ I(n) := \min_{A \subseteq \mathbb{Z}, |A| = n} \int^{1}_0 \left| \sum_{n \in A} e^{2\pi i n \theta}\right| \, d\theta \; ?\]
The problem of determining the exact minimum, known as the strong Littlewood conjecture,
is still unresolved (the conjecture is that the minimum occurs when $A$ is
an arithmetic progression). However the bound $I(n) \gg \log n$, which is tight up to a constant factor, was obtained independently by the second-named author \cite{konyagin1} and by McGehee, Pigno and Smith \cite{MPS} some 20 years ago.\vspace{11pt}

\noindent There are various natural variants of the problem that one might consider. In this paper we will be interested in the discrete setting of functions defined on the group $\Zp$.
If $f : \Zp \rightarrow \mathbb{R}$ is a function, then instead of discussing
exponential sums we speak of the Fourier transform of $f$. This is defined for $r \in \Zp$ by \[ \widehat{f}(r) := \E_{x \in \Zp}f(x) e(rx/p).\]
Here $e(\theta) = e^{2\pi i \theta}$, as is customary, and $\E$ denotes the averaging operator, which in this instance is the same as $p^{-1}\sum_{x \in \Zp}$. If $A \subseteq \Zp$ is a set then we write $1_A$ for its characteristic function. We will often be concerned with the Fourier transform $\widehat{1}_A$.

\begin{question}\label{q1} Let $p$ be a prime, and let $A \subseteq \Zp$ be a set of size $\lfloor p/2\rfloor$. What is $S(p)$, the minimum possible value of
\[ \Vert \widehat{1}_A \Vert_1 := \sum_{r \in \Zp} |\widehat{1}_A(r)|?\]
\end{question}
\noindent There is no serious reason for restricting attention to sets of size exactly $\lfloor p/2\rfloor$, though clearly if sets of any size are allowed then $S(p)$ can equal 1. An obvious lower bound for $S(p)$ is $p^{-1}\lfloor p/2\rfloor = 1/2 + O(p
^{-1})$, since this is the contribution to $\Vert \widehat{1_A} \Vert_1$ from the term $r = 0$. It is already non-trivial to obtain a bound of the form $S(p) \rightarrow \infty$. We will obtain a very weak bound of this sort in \S \ref{sec2}, followed in \S \ref{sec4} by the following stronger result.

\begin{theorem}\label{mainthm1} We have $S(p) \gg (\log p/\log \log p)^{1/3}$.\end{theorem}

\noindent The best upper bound we have comes by considering the set $A = \{1,\dots,\lfloor p/2\rfloor\}$, which illustrates that $S(p) \ll \log p$. By analogy with the strong Littlewood conjecture, one might guess that this represents the truth.\vspace{11pt}

\noindent Theorem \ref{mainthm1} will be obtained as a straightforward consequence of the next result, which is the main theorem in our paper.

\begin{theorem}\label{mainthm2} Let $f : \Zp \rightarrow \mathbb{R}$ be a function with $\E f = 0$ and $\Vert \widehat{f} \Vert_1 \leq 1$.
Then \[ \min_{x \in \Zp} |f(x)| = O\big((\log \log p/\log p)^{1/3}\big).\]
\end{theorem}
\noindent The condition $\Vert \widehat{f} \Vert_1 \leq 1$ may be thought of as expressing a (rather strong) kind of ``continuity'' or ``smoothness'' of $f$. Indeed one knows from harmonic analysis on the circle that if $f : \mathbb{R}/\Z \rightarrow \mathbb{R}$ is twice continuously differentiable then $\widehat{f}(m) \ll |m|^{-2}$ for $m \in \Z \setminus \{0\}$, and hence $\widehat{f}$ lies in $l^1(\Z)$. \vs

\ni The conclusion of Theorem \ref{mainthm2} may be thought of as a type of intermediate value theorem.\vs

\ni Taking $f(x) = cp/(p+1)\log p$ for $x = 0,1,\dots,\lfloor p/2\rfloor$ and $f(x) = -cp/(p-1)\log p$ for $x = -1,\dots,-\lfloor p/2\rfloor$, where $c$ is a suitably small constant, one sees that the bound in Theorem \ref{mainthm2} could not be improved
beyond $O(1/\log p)$. It may be that this is best possible.\vs

\noindent As a corollary of Theorem \ref{mainthm2}, we are able to obtain some information on a question of Gowers, the consideration of which was in fact the starting point of our investigations. Gowers asked whether, if $A \subseteq \Zp$ is a set of size $\lfloor p/2\rfloor$, there is some $x \in \Zp$ such that
\[ \left| |A \cap (A + x)| - p/4 \right| = o(p).\]
Note that $p/4$ is roughly the expected value of $|A \cap (A + x)|$. One interesting feature of this question, and of Theorem \ref{mainthm2}, is that similar statements do not hold in abelian groups $G$ more general than $\Zp$; for example if $G = \mathbb{F}_2^m$ and if $A$ is a hyperplane then clearly $|A \cap (A + x)|$ is always either $0$ or $|G|/2$. We answer a generalisation of Gowers' question in the affirmative, obtaining as a consequence of Theorem \ref{mainthm2} the following quantitative bound.

\begin{theorem}\label{mainthm3} Let $A,B \subseteq \Zp$ be subsets
of $\Zp$ with $|A| = \alpha p$, $|B| = \beta p$. Then there is some $x \in \Zp$ such that
\[ \big| |A \cap (B + x)| - \alpha \beta p \big| = O\left(p(\log \log p /\log p)^{1/3}\right).\]
\end{theorem}
\noindent If $f,g : \Zp \rightarrow \mathbb{R}$ are two functions then we define their convolution $f \ast g$ by
\[ f \ast g(x) := \E_y f(y) g(x - y).\] Note that
\[ |A \cap (B + x)| = p 1_A \ast 1_B^{\circ}(x),\] where we have written $f^{\circ}(x) = f(-x)$.
It is thus easy to see that Theorem \ref{mainthm3} is a special case of the following result:
\begin{theorem}\label{mainthm3.1} Let $f,g : \Zp \rightarrow \mathbb{R}$ be
two functions with $\Vert f \Vert_2, \Vert g \Vert_2\leq 1$. Then there is some $x \in \Zp$ such that
\[ \left| f \ast g(x)- \E f \E g\right|
 = O(\log \log p /\log p)^{1/3}.\]
\end{theorem}
\noindent For Gowers' original problem it is not clear how much one can hope for,
since there are no obvious examples of sets $A$ for which $|A \cap (A + x)|$
is always quite different to $p/4$.\vs

\ni In \S \ref{sec6} and in \S \ref{sec7} we prove the following results.
\begin{theorem}\label{mainthm3.2}
For any prime $p$ there is a function $f : \Zp \rightarrow \mathbb{R}$
such that $\E f=0$, $\Vert f \Vert_2 \leq 1$, and
\[ |f \ast f^{\circ}(x)| \geq c/\log p\] for all $x \in \Zp$.
\end{theorem}

\begin{theorem}\label{mainthm3.3}
For any prime $p\geq3$ there is a set $A \subseteq \Zp$,
$|A| = \lfloor p/2 \rfloor$, such that
\[ |A \cap (A + x) - p/4| \geq cp/\log p\log\log p\]
for all $x \in \Zp$.
\end{theorem}

\ni It would be of interest to close the gap between Theorem \ref{mainthm3} and \ref{mainthm3.3}, and also to obtain the correct bound in Theorem \ref{mainthm2}.

\section{The continuity argument, and a simple lower bound}\label{sec2}

\noindent The goal of this section is to motivate our main argument, and in so doing to prove a weak version of Theorem \ref{mainthm2}.
\begin{theorem}
Let $f : \Zp \rightarrow \mathbb{R}$ be a function with $\E f = 0$ and $\Vert \widehat{f} \Vert_1 \leq 1$.
Then \[ \min_{x \in \Zp} |f(x)| = o(1).\]
\end{theorem}
\ni By the Fourier inversion formula we have
\[ f(x) = \sum_{r \in \Zp} \widehat{f}(r) e(-rx/p).\] Since $\Vert \widehat{f} \Vert_1$ is so small, one might hope to have some rather direct control on this expression. Let us imagine, for a moment, that we know something even stronger about $\widehat{f
}$, namely that for some fairly small set $R \subseteq \Zp$ we have
\[ \sum_{r \notin R} |\widehat{f}(r)| \leq \epsilon/3.\]
By a well-known pigeonhole argument of Dirichlet, to be recalled in
\S \ref{sec3}, we may find some $t \neq 0$ such that the fractional parts
$\{tr/p\}$, $r \in R$, are all small (we define our fractional parts to lie between $-1/2$ and $1/2$). This implies that
$|e(tr/p) - 1|$ is always at most $\epsilon/3$, say. Thus we have
\begin{eqnarray*}
|f(x) - f(x+t)| & = & \sum_{r \in \Zp} \widehat{f}(r) \big( e(rx/p) - e(r(x+t)/p)\big) \\ & \leq & \sum_{r \in R} |\widehat{f}(r)| |1 - e(rt/p)| + 2 \sum_{r \notin R} |\widehat{f}(r)| \\ & \leq & \epsilon /3 + 2\epsilon/3 = \epsilon
\end{eqnarray*}
for all $x$. Since $\E f(x) = 0$, there must be some value of $x$ for which $f(x)$
and $f(x + t)$ have opposing signs. For this value of $x$ it is clear that $|f(x)| \leq \epsilon$.\vspace{11pt}

\ni\textit{Remark.} What we have done here is identify a length, $t$, along which $f$ is somewhat continuous in the sense that $|f(x + t) - f(x)|$ is always small.\vs

\noindent We now proceed to show how the above argument may be modified to give a weak version of Theorem \ref{mainthm2}. It need not be the case that most of the $L^1$-norm of $\widehat{f}$ is concentrated on a few modes, but by using a certain pigeonholing argument one may select a set $R$ of modes so that the contribution of $\sum_{r \notin R} \widehat{f}(r) e(-rx/p)$ is small for many $x$. \vspace{11pt}

\noindent Let $\epsilon > 0$, and define a sequence $(\delta_j)_{j=0}^{\infty}$
by $\delta_0 = 1$ and $\delta_{l+1} = \epsilon^{2 + 1/\delta_l}2^{-6/\delta_l - 4}$
for $l \geq 0$. Take $J = \lceil 3/\epsilon\rceil$. Then there exists $l \in \{0,\dots,J\}$ such that
\[ \sum_{r : \delta_{l+1} \leq |\widehat{f}(r)| < \delta_l } |\widehat{f}(r)| \leq \epsilon/3.\]
For this value of $l$, define
\[ f^{(1)}(x) =  \sum_{r : |\widehat{f}(r)| \geq \delta_l } \widehat{f}(r) e(-rx/p),\]
\[ f^{(2)}(x) =  \sum_{r : \delta_{l+1} \leq |\widehat{f}(r)| < \delta_l } \widehat{f}(r) e(-rx/p)\] and
\[ f^{(3)}(x) =  \sum_{r : |\widehat{f}(r)| < \delta_{l+1}} \widehat{f}(r) e(-rx/p).\]
By construction we have the estimate
\begin{equation}\label{eq2.1} \Vert f^{(2)} \Vert_{\infty} \leq \epsilon/3,\end{equation}
and also the bound
\begin{equation}\label{eq2.2} \Vert f^{(3)} \Vert_2^2 = \sum_{r : |\widehat{f}(r)| < \delta_{l+1}} |\widehat{f}(r)|^2 \leq \delta_{l+1}.\end{equation}
Write $R$ for the set of frequencies occurring in $f^{(1)}$, that is to say
$R := \{ r : |\widehat{f}(r)| \geq \delta_l \}$, and set $k := |R|$. By the assumption on $\Vert \widehat{f} \Vert_1$ we have $k \leq 1/\delta_l$. Consider the Bohr set
\[ B(R,\eta) := \{x \in \Zp : |\{ x r /p \}|
\leq \eta \; \mbox{for all $r \in R$}\}.\]
By Dirichlet's argument we have the bound
\begin{equation}\label{eq2.3} |B(R,\eta)| \geq p\eta^k \geq p\eta^{1/\delta_l}.\end{equation}
For details, and more properties of Bohr sets, see \S \ref{sec3} (particularly Lemma \ref{lem4.00} (i) and its proof). \vspace{11pt}

\noindent Let $B := B(R,\epsilon/48)$. Suppose that $B\setminus\{0\}$ is nonempty,
and select an element $t$ from it. By the same argument we had earlier, we have
\[
|f^{(1)}(x) - f^{(1)}(x+t)| \leq  \sum_{r \in R} |\widehat{f^{(1)}}(r)| |1 - e(rt/p)|  \leq \epsilon/6.\]
Once again, since $\E f^{(1)}(x) = 0$, there is $x_0$ such that $f^{(1)}(x_0)$ and $f^{(1)}(x_0 + t)$ have opposing signs, which implies that $|f^{(1)}(x_0)| \leq \epsilon/6$. For any $t \in B$ we then have $|f^{(1)}(x_0 + t)| \leq \epsilon/3$. \vspace{11
pt}

\noindent Set $X := x_0 + B$. Combining the above observations with \eqref{eq2.1} we have, if $x \in X$, the inequality
\[ |f(x)| \leq 2\epsilon/3 + |f^{(3)}(x)|.\]
To conclude the argument it remains to show that there is some $x \in X$ for which $|f^{(3)}(x)| \leq \epsilon/3$. To this end, note that if this were not the case then we should have
\begin{equation}\label{eq2.4} \Vert f^{(3)} \Vert_2^2 = \E_{x \in \Zp} f^{(3)}(x)^2  \geq \frac{|B|}{p}\E_{x \in X}f^{(3)}(x)^2 > |B|\epsilon^2/9 \geq \epsilon^{2 + 1/\delta_l}2^{-6/\delta_l - 4}p,\end{equation}
which contradicts \eqref{eq2.2}.\vspace{11pt}

\noindent The only condition that must be satisfied to make this argument work
is that $B\setminus\{0\}$ must be non-empty. By dint of \eqref{eq2.3} this will be so if $p(\epsilon/48)^{1/\delta_J} > 1$. Since $\delta_J$ is subject to a lower bound
depending only on $\epsilon$, this can certainly be satisfied with
$\epsilon = o(1)$, and in fact one can take $3/\epsilon = \log_{*}(p) + O(1)$. This is a very slowly growing function, of course: the reason for this is that $1/\delta_J$ looks like a tower of exponentials of height about $J$, this being forced by the need to have $\delta_{j+1}$ exponentially smaller than $\delta_j$.\vspace{11pt}

\noindent Before moving on, let us highlight the weakness in the above argument which we will be able to tighten up in subsequent sections, leading to a superior bound. It lies in the fact that \eqref{eq2.4} represents not just a lower bound for $\Vert f^{(3)} \Vert_2^2$, but in fact for $\E_{x \in X} |f^{(3)}(x)|^2$. In order to exploit this properly one needs to modify \eqref{eq2.2} appropriately. Writing $\beta(x) := 1_B(x)/|B|$ for the \textit{normalised Bohr cutoff} associated to
 $B$, our interest lies in finding a Fourier expression for
\[ \E_{x \in \Zp} |f^{(3)}(x)|^2 \beta(x + x_0) .\] This is not hard; indeed we have
\begin{eqnarray}\nonumber
\E_x |f^{(3)}(x)|^2 \beta(x + x_0) & = & \E_x \sum_{r,r'} \widehat{f^{(3)}}(r)\overline{\widehat{f^{(3)}}(r')} e((r' - r)x/p) \beta(x + x_0) \\ & = & \sum_{r,r'}
\widehat{f^{(3)}}(r)\overline{\widehat{f^{(3)}}(r')} \widehat{\beta}(r - r') e((r - r')x_0/p), \label{eq2.6}
\end{eqnarray}
an expression we will find very useful later on. In exploiting it we will require a few facts concerning Bohr sets and their Fourier transforms. We develop these in the next section.

\section{Properties of Bohr sets}
\label{sec3}
\noindent The material of this section is somewhat technical, but is rapidly becoming part of the foundation of additive combinatorics. It has its origins in the work of Bourgain \cite{Bou} on 3-term arithmetic progressions. See \cite{tao-exposition} for
a recent exposition of this work. For examples of subsequent work in a similar vein, see either \cite{green1} or \cite{tao1}, and for a more leisurely discussion of Bohr sets and the need for the results of this section, see \cite{green-finite-fields}. We
 work in the context of an arbitrary finite abelian group $G$, since this is most natural, though our interest is in the case $G = \mathbb{Z}/p\mathbb{Z}$. In fact, a lot of the material in this paper is valid in this more general setting, the key exception being the ``discrete intermediate value theorem'' of \S \ref{sec1}. The failure of this, of course, means that the key theorems as stated in \S \ref{sec1} are only valid when $G = \mathbb{Z}/p\mathbb{Z}$.\vs

\ni If $f : G \rightarrow \mathbb{R}$ is a function and if $\gamma \in G^*$ is a character then we define the Fourier transform $\widehat{f}$ by
\[ \widehat{f}(\gamma) := \E_{x \in G} f(x) \gamma(x).\]
It is easy to tie this in with the definition we had in the case $G = \Zp$; in that group, all characters are of the form $x \mapsto e(rx/p)$ for some $r \in \Zp$.\vspace{11pt}

\noindent Now if $\Gamma \subseteq G^{\ast}$ is a set of $d$ characters, and if $\epsilon > 0$, then we define the Bohr set $B = B(\Gamma,\epsilon)$ by
\[ B(\Gamma,\epsilon) := \{ x \in G : |\textstyle \frac{1}{2\pi}\displaystyle \arg\gamma(x)| \leq \epsilon \;\; \mbox{for all $\gamma \in \Gamma$}\}.\]We also define $\beta = \beta_{\Gamma,\epsilon}$ (the \emph{normalised Bohr cutoff} associated to $B$) by $\beta(x) := 1_{B}(x)/|B|$.\vs

\ni The next lemma details some well-known facts concerning the size of Bohr sets.

\begin{lemma}\label{lem4.00} Let $\Gamma \subseteq G^*$ be a set of $d$ characters, and let $\delta > 0$. Then we have the bounds
\begin{enumerate}
\item[\textup{(i)}] $|B(\Gamma,\delta)| \geq \delta^d |G|$;
\item[\textup{(ii)}] $|B(\Gamma,2 \delta)| \leq 5^d|B(\Gamma,\delta)|$.
\end{enumerate}
\end{lemma}
\proof For any $\eta$, write $S_{\eta}$ for the set of all
$y = (y_1,\dots,y_d) \in \mathbb{R}^d/\mathbb{Z}^d$ for which
$\Vert y_j\Vert_{\infty} \leq \eta/2$ for all $j = 1,\dots,d$ (to define the $\Vert \cdot \Vert_{\infty}$ norm on $\mathbb{R}/\Z$ we identify it with $(-\frac{1}{2},\frac{1}{2}]$).
Let $\Gamma=\{(\gamma_1,\dots,\gamma_d)\}$. Now if $x \in G$ write \[ v(x) \; = \; \textstyle \frac{1
}{2\pi} \displaystyle(\arg \gamma_1(x),\dots,\arg \gamma_d(x)) \; \in \;\mathbb{R}^d/\mathbb{Z}^d.\]
If $v(x)$ and $v(x')$ both lie in some translate $a + S_{\delta}$ then $x - x' \in B(\Gamma,\delta)$,
and so for fixed $x' \in G\cap v^{-1}((a + S_{\delta}))$ the map $x \mapsto x - x'$ defines an injection from $G \cap v^{-1}(a + S_{\delta})$ to $B(\Gamma,\delta)$. Hence
\begin{equation}\label{eq476} |v(G) \cap (a + S_{\delta})| \; \leq \; |B(\Gamma,\delta)|.\end{equation}
\textit{Proof of} (i). By a simple averaging there is some translate $a + S_{\delta}$ such that
\[ |v(G) \cap (a + S_{\delta})| \; \geq \; |S_{\delta}||v(G)| \; = \; |S_{\delta}| \; = \; \delta^d|G|.\] The result is now immediate from \eqref{eq476}.\vs

\ni\textit{Proof of} (ii). From \eqref{eq476} one has
\[ \left|v(B(\Gamma,2 \delta)) \cap (a + S_{\delta})\right| \; \leq \;\left|v(G) \cap (a + S_{\delta})\right| \; \leq \; |B(\Gamma,\delta)|.\]
Now $v(B(\Gamma,2\delta)) \cap (a + S_{\delta})$ is empty unless $a \in S_{5\delta}$, and so
\begin{equation}\boxeq
|B(\Gamma,2\delta)| \; = \; \frac{1}{|S_{\delta}|} \int_{\mathbb{R}^d/\mathbb{Z}^d} \left|v(B(\Gamma,2 \delta)) \cap (a + S_{\delta})\right| da \; \leq \; |B(\Gamma,\delta)| \cdot \frac{|S_{5\delta}|}{|S_{\delta}|} \; = \; 5^d \cdot |B(\Gamma,\delta)|.
\end{equation}

\noindent Suppose now that $x \in B(\Gamma,\epsilon)$ and that $x' \in B(\Gamma,\epsilon')$, where $\epsilon' \ll \epsilon$. Then $x + x' \in B(\Gamma,\epsilon + \epsilon')$. Assuming that $B(\Gamma,\epsilon + \epsilon')$ is not much bigger than $B(\Gamma
,\epsilon)$, this means that $B(\Gamma,\epsilon) + B(\Gamma,\epsilon') \approx B(\Gamma,\epsilon)$. After passing to the associated normalised cutoff functions $\beta$ and $\beta'$, we might anticipate that $\beta \ast \beta' \approx \beta$. When true this property, which says that Bohr cutoffs are roughly invariant under convolution by narrower Bohr cutoffs, is extremely useful.\vspace{11pt}

\noindent Unfortunately the property does not always hold, since the size of $B(\Gamma,t)$ need not vary very smoothly with $t$. For example, suppose that $G = \mathbb{F}_5^n$ and that $\delta_1 < 1/5 < \delta_2$. Then if $|\Gamma|$ consists of $d$ linearly independent characters we have $|B(\Gamma, \delta_1)| = 5^{-d}|G|$ whilst $|B(\Gamma,\delta_2)| = (3/5)^d |G|$. \vspace{11pt}

\noindent Bourgain showed, and we shall repeat his argument, that for a fixed $\Gamma$ there are a plentiful supply of $\epsilon$ such that the function $|B(\Gamma,t)|$ does vary quite regularly for $t \approx \epsilon$, in the following sense.

\begin{definition}[Regular value]
Let $\Gamma$ be a set of $d$ characters. We say that $B(\Gamma,\epsilon)$ is \textup{regular value}, or that $\epsilon$ is a \textup{regular value} for $\Gamma$, if
\[ 1 - 100d |\kappa| \leq \frac{|B(\Gamma,(1 + \kappa)\epsilon)|}{|B(\Gamma,\epsilon)|}  \leq  1 + 100d|\kappa| \] whenever $|\kappa| \leq 1/100d$.
\end{definition}

\begin{proposition}\label{reg-prop} Let $\Gamma$, $|\Gamma| = d$, be a fixed set of characters. Let $\delta \in (0,1)$. Then there is $\epsilon \in [\delta,2\delta)$ which is regular for $\Gamma$.
\end{proposition}
\proof
Set $f(\alpha) = |B(\Gamma,2\alpha\delta)|$. Then
$f(\alpha)$ is a non-decreasing function on $[\frac{1}{2},1]$ and, by Lemma \ref{lem4.00},
\begin{equation} \label{fineq} f(1)  \leq  5^d f(1/2).\end{equation}  We wish to show that there is $\alpha \in [1/2,1]$ such that $1 - 100d |\kappa|  \leq  f((1 + \kappa)\alpha)/f(\alpha)  \leq  1 + 100 d |\kappa|$ for all $|\kappa| \leq 1/100d$. Suppose
 then that this is false. Observe that $\frac{1}{1 - x} \geq 1 + x$ when $x \geq 0$; hence for every $\alpha \in [1/2,1]$ there is $t_{\alpha} \in [0,1/100d]$ such that
\begin{equation}
\label{eq2} \left| \frac{f \left( (1 + t_{\alpha}) \alpha \right)}{f \left( (1 - t_{\alpha}) \alpha \right)} \right|  \geq  1 + 100 d t_{\alpha}
 \geq  e^{50 d t_{\alpha}}, \end{equation}
the last step following because $1 + x \geq e^{x/2}$ for $x \leq 1$.\vspace{11pt}

\noindent The next lemma is related to the Vitali covering lemma, but is so simple that we provide an independent proof. For an anecdotal discussion see \cite{croft}.
\begin{lemma}
\label{lemma6}
Suppose a finite collection of closed intervals $I_1,\ldots,I_k$ covers $[0,1]$. Then we can pick a subcollection $I_{i_1},\ldots,I_{i_m}$ whose members are disjoint except possibly at their endpoints, with total measure at least $1/2$.
\end{lemma}
\proof Without loss of generality suppose that the collection $I_1,\ldots,I_k$ is minimal in that if any $I_j$ is removed, the intervals no longer cover $[0,1]$. It is then easy to see that no point $x$ lies in three of the $I_j$, because there are two in
tervals $I_r$ and $I_s$ containing $x$ such that any other $I_t$ containing $x$ lies in $I_r \cup I_s$. But it is then easy to describe the form of the intervals exactly. Suppose that the $I_j = [a_j,b_j]$ with $a_1 \leq a_2 \leq \ldots \leq a_k$. Then \[
 a_1 \leq a_2 \leq b_1 \leq a_3 \leq b_2 \leq a_4 \leq \ldots \leq b_{k-1} \leq b_k. \]
It follows that the two collections $I_1 \cup I_3 \cup \ldots$ and $I_2 \cup I_4 \cup \ldots$ contain disjoint intervals. The result is now obvious. \endproof\vspace{11pt}

\noindent To apply Lemma \ref{lemma6}, recall \eqref{eq2}. By compactness we may take a finite set \[ \{\alpha_1,\ldots,\alpha_k\} \subseteq [ \textstyle\frac{1}{2} + \frac{1}{100d},1 - \frac{1}{100d} ] \] such that the intervals
$\left[ \left(1 - t_{\alpha_i}\right)\alpha_i,\left(1 + t_{\alpha_i}\right)\alpha_i \right]$ cover $\left[ \frac{1}{2} + \frac{1}{100d},1  -  \frac{1}{100d}  \right]$. Since  $t_{\alpha_i}  \leq  1/100d$, all of these intervals are contained in
$[\frac{1}{2},1]$. By Lemma \ref{lemma6},
 we can pick a disjoint subcollection of measure at least $\frac{1}{4}\left(1 - 1/100d\right) > \frac{1}{5}$. Letting these intervals correspond to $\{ \alpha_{1},\ldots,\alpha_{l} \}$, one has
$2 \sum_{i = 1}^l \alpha_i t_{\alpha_i}  >  1/5$,
and so $\sum_{i = 1}^l t_{\alpha_i}  >  1/10$.
Using this in \eqref{eq2} gives
\[ \prod_{i = 1}^l \left| \frac{f \left( (1 + t_{\alpha_i}) \alpha_i \right)}{f \left( (1 - t_{\alpha_i}) \alpha_i \right)} \right| \; \geq \; e^{50d \sum_{i} t_{\alpha_i}} \; > \; e^{5d}.\]
However the left hand side is at most $f(1)/f(1/2)$, and hence by \eqref{fineq} we have
\[ 5^d \; > \; e^{5d}.\]
This is a contradiction, and Proposition \ref{reg-prop} is established. \endproof\vspace{11pt}

\ni\textit{Remark.} Actually, Proposition \ref{reg-prop} easily follows from
weak type estimates for the Hardy--Littlewood maximal function. Let
$\nu$ be a bounded non-decreasing function $\mathbb{R} \to \mathbb{R}$
with an associated Radon---Stieltjes measure
\[\nu[a,b)=\nu(b)-\nu(a)\quad(a<b).\]
Note that
\[\nu(\mathbb{R})=\lim_{b\to\infty}\nu[-b,b).\]
The Hardy--Littlewood maximal function is defined as
\[M_\nu(t)=\sup_{a<t<b}\frac{\nu(b)-\nu(a)}{b-a}.\]
The weak type estimates for $M_n(u)$ claim that for any $\alpha>0$
\begin{equation}\label{wti}
\mu\{t:\,M_\nu(t)>\alpha\}\le\frac{2\nu(\mathbb{R})}{\alpha}
\end{equation}
where $\mu$ is the Lebesgue measure. For an absolutely continuous $\nu$
this inequality is classical; a general case was established in \cite{luk}.
One can prove \ref{reg-prop}
applying (\ref{wti}) to the function $\nu$ such that $\nu(x)=\log f(e^x)$
for $x\in[-\log 2, 0]$, $\nu(x)=\nu(-\log 2)$ for $x<-\log 2$, and
$\nu(x)=\nu(0)$ for $x>0$. \vs

\begin{lemma}\label{lem6.7} Let $\Gamma \subseteq G^*$ be a set of $d$ characters, and let $\delta,\epsilon > 0$. Suppose that $\eta$ is a regular value for $\Gamma$, and let $\beta = \beta_{\Gamma,\eta}$ be the normalised Bohr cutoff associated to
$B(\Gamma,\eta)$. Suppose that $y \in B(\Gamma,\eta')$, where $\eta' \leq \epsilon\eta/200d$. Then
\[ \E_{x \in G} |\beta(x + y) - \beta(x)|  \leq \epsilon.\]
\end{lemma}
\proof Set $\delta := \epsilon/200d$. If $\beta(x+y) - \beta(x) \neq 0$ then
$x \in B(\Gamma, \eta(1 + \delta)) \setminus B(\Gamma,\eta(1 - \delta))$. Thus, since $B = B(\Gamma,\eta)$ is regular, we have
\[ \sum_{x} |1_B(x+y) - 1_B(x)| \leq 200d\delta |B|,\]
which immediately implies the result.\endproof\vs

\ni The following corollary will be very useful. It tells us that if $\widehat{\beta}(\gamma)$ is moderately large, then $\gamma(y) \approx 1$ for $y$ in a smaller Bohr set $B(\Gamma,\eta')$.

\begin{lemma}\label{lem3.6}
Let $B(\Gamma,\eta)$ be a regular Bohr set with normalised cutoff $\beta = \beta_{\Gamma,\eta}$. Suppose that $\kappa_1,\kappa_2 > 0$. Suppose that $\gamma$ is a character for which  $|\widehat{\beta}(\gamma)| \geq \kappa_1$, and that $y \in B(\Gamma,\eta
')$, where $\eta' \leq \kappa_1 \kappa_2 \eta/200d$. Then $|1 - \gamma(y)| \leq \kappa_2$.
\end{lemma}
\proof We have
\[ \kappa_1 |1 - \gamma(y)| \leq |\widehat{\beta}(\gamma)||1 - \overline{\gamma(y)}| = \big|\E_{x \in G} (\beta(x+y) - \beta(x))\gamma(x)\big|,\]
which is at most $\kappa_1 \kappa_2$ by Lemma \ref{lem6.7}. The result follows. \endproof\vspace{11pt}

\ni\textit{Remark.} By employing appropriate smoothing devices, one may obtain a reasonably good quantitative description of the set of $\gamma$ for which $|\widehat{\beta_{\Gamma,\epsilon}}(\gamma)| \geq \delta$ when $B(\Gamma,\epsilon)$ is regular. If $
\Gamma = \{\gamma_1,\dots,\gamma_d\}$ these characters will all be of the form $\gamma_1^{a_1}\dots \gamma_d^{a_d}$, where the $a_i$ are integers bounded in absolute value by some $F(\delta)$. We will not need to use such information in this paper.

\section{Proofs of the main theorems}\label{sec4}
\ni Our first task in this section will be to supply a proof of Theorem \ref{mainthm2}. We begin with some motivating remarks which, it will turn out, take us some distance into the proof itself. Let $f : \Zp \rightarrow \mathbb{R}$ be a function with
$\Vert \widehat{f} \Vert_1 \leq 1$, as in the statement of that theorem. Suppose to begin with that we have defined a splitting $f = f^{(1)} + f^{(2)}$ where
\[ f^{(1)}(x) := \sum_{\gamma \in \Gamma} \widehat{f}(\gamma) \overline{\gamma(x)}\]
for some reasonably small set $\Gamma$. If $f^{(1)}(x_0)$ is small then $f^{(1)}(x + x_0)$ will also be small provided $x \in B(\Gamma,\eta)$, for suitably small $\eta$. Our first observation is that much the same would be true if
$\Gamma = \{\gamma_1,\dots,\gamma_d\}$ were replaced by some set
\[ \{ \gamma_1^{a_1}\dots \gamma_d^{a_d} : |a_i| \leq M\},\] consisting of small combinations of elements of $\Gamma$. By remarks made at the end of \S \ref{sec3}, the characteristic function of this set rather resembles the function
$\widehat{\beta_{\Gamma,\eta}}$, for suitable $\eta$. Once this is noticed, one might think to redefine
\[ f^{(1)} := \sum_{\gamma} \widehat{\beta_{\Gamma,\eta}}(\gamma) \widehat{f}(\gamma) \overline{\gamma(x)} = f \ast \beta_{\Gamma,\eta},\]
which looks extremely natural. To begin one takes $\Gamma = \emptyset$; we add frequencies according to an iterative procedure.
\begin{lemma}\label{lem4.1}
Suppose that $B(\Gamma,\eta)$ is regular, that $|f^{(1)}(x_0)| \leq \epsilon$,
and that $x \in B(\Gamma,\eta')$ for some $\eta' \leq \epsilon^2 \eta/200d$.
Then $|f^{(1)}(x + x_0)| \leq 3\epsilon$.
\end{lemma}
\proof We have
\[
|f^{(1)}(x + x_0) - f^{(1)}(x_0)|  =   \big|\sum_\gamma \widehat{f}(\gamma)\widehat{\beta_{\Gamma,\eta}}(\gamma)\gamma(x_0)(\gamma(x) - 1)\big| \leq  \sup_{\gamma} \big|\widehat{\beta_{\Gamma,\eta}}(\gamma)\big| \big|1 - \gamma(x)\big|.\]
Now if $\big|\widehat{\beta_{\Gamma,\eta}}(\gamma)\big| \leq \epsilon$ then clearly
\begin{equation}\label{eq3.87} \big|\widehat{\beta_{\Gamma,\eta}}(\gamma)\big| \big|1 - \gamma(x)\big| \leq 2\epsilon.\end{equation}
If, however, $\big|\widehat{\beta_{\Gamma,\eta}}(\gamma)\big| \geq \epsilon$ and if $x \in B(\Gamma,\eta')$ then Lemma \ref
{lem3.6} tells us that $|1 - \gamma(x)| \leq \epsilon$, and so \eqref{eq3.87} holds in this case too.\endproof\vspace{11pt}

\ni Now define
\begin{equation}\label{eq4.4a} f^{(2)} := f - f^{(1)}.\end{equation} Suppose that $f^{(1)}(x_0)$ is small. By Lemma \ref{lem4.1} we know that $f^{(1)}(x + x_0)$ is small for all $x \in B(R,\eta')$, provided $\eta'$ is not too large. Our hope is that $f^{(1)}(x + x_0)$ will also be small for some such $x$. To this end we look at
\begin{equation}\label{eq4.3} E := \E_{x \in G}  |f^{(2)}(x + x_0)|^2 \beta_{\Gamma,\eta'}(x) \leq \sum_{\gamma,\gamma'} |\widehat{f^{(2)}}(\gamma)||\widehat{f^{(2)}}(\gamma')||\widehat{\beta_{\Gamma,\eta}}(\gamma - \gamma')|,
\end{equation}
the inequality being an immediate consequence of \eqref{eq2.6}. Write
\[ U := \{ \chi : |\widehat{\beta_{\Gamma,\eta'}}(\chi)| \geq \epsilon^2\}.\] Observing that
\[ \sum_{\gamma} |\widehat{f^{(2)}}(\gamma)| = \sum_{\gamma} |\widehat{f}(\gamma)||1 - \widehat{\beta_{\Gamma,\eta}}(\gamma)| \leq \Vert \widehat{f} \Vert_1\sup_{\gamma} | 1 - \widehat{\beta_{\Gamma,\eta}}(\gamma)| \leq 2,\]
it is immediate from \eqref{eq4.3} that
\begin{equation}\label{eq4.4} E \leq 4\epsilon^2  + \sum_{\substack{\gamma,\gamma' \\ \gamma - \gamma' \in U}} |\widehat{f^{(2)}}(\gamma)||\widehat{f^{(2)}}(\gamma')||\widehat{\beta_{\Gamma,\eta'}}(\gamma - \gamma')| \leq 4\epsilon^2  + 2\sup_{\gamma}
\sum_{\gamma' : \gamma' - \gamma \in U} |\widehat{f^{(2)}}(\gamma')|.\end{equation}
Now if
\begin{equation}\label{eq4.5b} \sup_{\gamma} \sum_{\gamma' : \gamma' - \gamma \in U} |\widehat{f^{(2)}}(\gamma')| \leq \epsilon^2\end{equation}
then $E \leq 6\epsilon^2$, and there is indeed $x \in B(\Gamma,\eta')$ such that
$|f^{(2)}(x+x_0)| \leq \epsilon\sqrt{6} \leq 3 \epsilon$. For this $x$, we have
$|f(x+x_0)| \leq |f^{(1)}(x+x_0)| + |f^{(2)}(x+x_0)| \leq 6\epsilon$. \vspace{11pt}

\noindent If, on the other hand, there is an $\gamma$ such that
\begin{equation}\label{eq4.5a} \sum_{\gamma' : \gamma' - \gamma \in U} |\widehat{f^{(2)}}(\gamma')| \geq \epsilon^2\end{equation} then we must find an alternative argument. It is tempting to somehow ``add'' the offending characters $u + \gamma$, $u \in U$
, into the definition of $f^{(1)}$. Thus we might consider a new set $\widetilde{\Gamma} := \Gamma \cup \{\gamma\}$ and some $\widetilde{\eta} \ll \eta$ such that $\widehat{\beta_{\widetilde{\Gamma},\widetilde{\eta}}}(s) \approx 1$ for $s = u + \gamma$.
\vspace{11pt}

\noindent Once this is done, one may iterate the whole procedure starting with \eqref{eq4.4a}. If things are done carefully, the obstacle \eqref{eq4.5a} cannot be encountered more than $\epsilon^{-2}$ times since each new instance corresponds to at least
$\epsilon^2$ of the total $L^1$-mass of $\widehat{f}$, which we are assuming is at most $1$. Care must be taken to ensure that these portions of $L^1$-mass are all disjoint, but we will see later on that this can be arranged.\vspace{11pt}

\noindent Working out the details of the above argument, one obtains a bound $\min_x |f(x)| = O((\log p)^{-1/4 + \varepsilon})$ in Theorem \ref{mainthm2}. To obtain the superior exponent $\frac{1}{3}$, we exploit the fact that it is not necessary for $f^{(2)}$ itself to satisfy \eqref{eq4.5b}; we may perturb $f^{(2)}$ by an arbitrary function $f^{(3)}$ for which $\Vert \widehat{f^{(3)}} \Vert_1 \leq \epsilon$, since then
\begin{equation}\label{eq4.12} |f^{(3)}(x)| = \big| \sum_{\gamma} \widehat{f^{(3)}}(\gamma) \overline{\gamma(x)} \big| \leq \Vert \widehat{f^{(3)}} \Vert_1 \leq \epsilon\end{equation} for an arbitrary $x$.
Thus in effect we wish to write \begin{equation}\label{eq4.11} f = f^{(1)} + f^{(2)} + f^{(3)},\end{equation} where $f^{(2)}$ satisfies \eqref{eq4.5b} and \begin{equation}\label{eq4.8}\Vert \widehat{f^{(3)}} \Vert_1 \leq \epsilon.\end{equation} This, it transpires, can be achieved by a slight variant of the above iterative scheme in which the number of iterations is reduced to just $\epsilon^{-1}$. \vs

\ni We now turn to the proof proper of Theorem \ref{mainthm2}, beginning with a lemma summarising the above discussion.
\begin{lemma}\label{lem4.14}
Let $B = B(\Gamma,\eta)$ and $B' = B(\Gamma,\eta')$ be regular Bohr sets with normalised Bohr cutoffs $\beta= \beta_{\Gamma,\eta}$ and $\beta' = \beta_{\Gamma,\eta'}$. Suppose that $\eta' \leq \epsilon^2 \eta/200d$. Set \[ U := \{ \chi : |\widehat{\beta'}
(\chi)| \geq \epsilon^2\}.\]
Suppose that $f = f^{(1)} + f^{(2)} + f^{(3)}$, where $f^{(1)} = f \ast \beta$, that
\[ \sup_{\gamma} \sum_{\gamma' : \gamma' - \gamma \in U} |\widehat{f^{(2)}}(\gamma')| \leq \epsilon^2 \] and that
$\Vert \widehat{f^{(3)}} \Vert_1 \leq \epsilon$. Suppose that $\min_x |f^{(1)}(x)| \leq 2\epsilon$.
Then $\min |f(x)| \leq 8\epsilon$.
\end{lemma}
\proof Suppose that $|f^{(1)}(x_0)| \leq 2\epsilon$. Then Lemma \ref{lem4.1}
tells us that for any $x \in B(R,\eta')$ we have $|f^{(1)}(x + x_0)| \leq 4\epsilon$. Now by the discussion leading to \eqref{eq4.5b} we know that there is $x \in B(R,\eta')$ for
which $|f^{(2)}(x + x_0)| \leq 3\epsilon$, and finally from \eqref{eq4.12} we have $|f^{(3)}(x + x_0)| \leq \epsilon$ for any $x$ whatsoever. The result follows. \endproof\vspace{11pt}

\noindent The next lemma asserts that there is a decomposition of the type discussed in Lemma \ref{lem4.14}.
\begin{lemma}\label{lem4.33}
Let $\epsilon > 0$. Then there is a set $\Gamma$, $|\Gamma| \leq 2/\epsilon^2$, of characters, an $\eta \geq (\epsilon^9/2^{21})^{2/\epsilon}$ and an $\eta' \geq 2^{-10}\epsilon^3 \eta$ such that the hypotheses of Lemma \ref{lem4.14} are satisfied, with the possible exception of the inequality $\min_x |f^{(1)}(x)| \leq \epsilon$.
\end{lemma}
\proof We will have a kind of double iterative scheme, the \texttt{STEP}s of which will be indexed by pairs $(j,i)$. The two coordinates of $\stp$ will be denoted $\stp_1$ and $\stp_2$. The algorithm will proceed by first incrementing $\stp_2$, that is to
 say according to the following scheme:
\[ (0,0) \rightarrow (0,1) \rightarrow \dots \rightarrow (0,l_0) \rightarrow (1,0) \rightarrow (1,1) \rightarrow \dots \rightarrow (1,l_1) \rightarrow (2,0) \rightarrow \dots\]
Whilst $\stp_1 = j$ there will be a set $\Gamma_{j}$ of characters and regular
Bohr sets $B(\Gamma_{j},\eta_{j})$, $B(\Gamma_{j},\eta'_{j})$ with normalised Bohr cutoffs $\beta_{j} := \beta_{\Gamma_{j},\eta_{j}}$ and $\beta'_j := \beta_{\Gamma_{j},\eta'_{j
}}$. The parameters $\eta_{j}$ and $\eta'_j$ will be related by the inequalities
\begin{equation}\label{eq4.20} \epsilon^2 \eta_{j}/400d \leq \eta'_j \leq \epsilon^2 \eta_{j}/200d,
\end{equation}
in order that we may apply Lemma \ref{lem4.14}. In fact $\eta'_{j}$ can simply
be an arbitrary parameter in this range for which $B(\Gamma_j,\eta'_{j})$ is regular. There will also be disjoint sets $S_{j,i}$ of characters, $i = 1,2,\dots$. When $\stp = (j,
i)$ we will have occasion to discuss the union
\[ \Omega_{j} := \bigcup_{j' < j} \bigcup_{i\leq l_{j'}}S_{j',i}\]
of all characters which were defined whilst $\stp_1 < j$. Things will have been arranged so that $B(\Gamma_{j},\eta_{j}) \subseteq B(\Omega_{j},\epsilon)$; the set $\Omega_{j}$ should be thought of as containing small combinations of elements of $\Gamma_{j}$, and in fact the proof could be set up in such a way that this is made explicit (see the remarks at the end of \S \ref{sec3}). In proceeding from $\stp = (j,0)$ to $\stp = (j,i)$ we will have accrued a collection $\{\gamma_{j,1},\dots,\gamma_{j,i}\}$
of characters; these will be used to define $\Gamma_{j+1}$.\vspace{11pt}

\noindent To initialize the iteration set $\stp = (0,0)$, and define $\Gamma_0 = \emptyset$ and $\eta_0 = 1$.  Now suppose that $\stp = (j,i)$. Define
\[ f^{(1)}_j := f \ast \beta_{j},\] and set $g_j := f - f^{(1)}_j$. \vspace{11pt}

\ni Define $U_{j} := \{ s : |\widehat{\beta'_j}(s)| \geq \epsilon^2\}$.
We will have reached $\stp = (j,i)$ by incrementing $\stp_2$, starting from $\stp = (j,0)$. During this process we will have defined disjoint sets $S_{j,1},S_{j,2},\dots,S_{j,i}$ of frequencies, and also a collection $\{\gamma_{j,1},\dots,\gamma_{j,i}\}$.
 The construction will be such that
\begin{equation}\label{eq4.16}
S_{j,i} \subseteq \gamma_{j,i} + U_{j}
\end{equation}
for all $i$.
We begin by asking whether or not
\begin{equation}\label{eq4.7} \sum_{\gamma' \in S_{j,1} \cup S_{j,2} \cup \dots \cup S_{j,i}} |\widehat{g_j}(\gamma')| \geq \epsilon.\end{equation}
If so, it is time to increment $\stp_1$ and to define $\Gamma_{j+1}$ and $\eta_{j+1}$. We set $\stp = (j+1,0)$, $\Gamma_{j+1} := \Gamma_j \cup \{\gamma_{j,1},\dots,\gamma_{j,i}\}$ and choose $\eta_{j+1}$,
\[  \epsilon^3\eta'_j/800|\Gamma_j| \leq \eta_{j+1} \leq \epsilon^3\eta'_j/400|\Gamma_j|,\] in such a way that the Bohr set $B(\Gamma_{j+1},\eta_{j+1})$ is regular. Clearly $\eta_{j+1} \leq \min(\eta_{j},\epsilon/2)$. Note that
\[ \Omega_{j+1} = \Omega_{j} \cup S_{j,1} \cup S_{j,2} \cup \dots \cup S_{j,i}.\]
\noindent\textit{Claim 1.} We have the inclusion \[B(\Gamma_{j+1},\eta_{j+1}) \subseteq B(\Omega_{j+1},\epsilon).\]
\proof Clearly \[ B(\Gamma_{j+1},\eta_{j+1}) \subseteq B(\Gamma_j,\eta_{j}) \subseteq B(\Omega_{j},\epsilon) \subseteq B(\Omega_{j+1},\epsilon),\] and so it suffices to check that $B(\Gamma_{j+1},\eta_{j+1}) \subseteq B(S_{j,i},\epsilon)$. For this, in view of \eqref{eq4.16}, we only need confirm that \[ B(\Gamma_{j+1},\eta_{j+1}) \subseteq B(\gamma_{j,i},\epsilon/2) \cap B(U_{j},\epsilon/2)\] for each $i$. The inclusion $B(\Gamma_{j+1},\eta_{j+1}) \subseteq B(\gamma_{j,i},\epsilon/2)$ is immediate from t
he fact that $\eta_{j+1} \leq \epsilon/2$. To see that $B(\Gamma_{j+1},\eta_{j+1}) \subseteq B(U_{j},\epsilon/2)$, we use Lemma \ref{lem3.6}. Indeed if $\gamma \in U_{j}$ and if $y \in B(\Gamma_j,\eta_{j+1})$ then by that lemma we have $|1 - \gamma(y)| \leq \epsilon/2$, which is exactly what we need.\vs

\noindent Now if \eqref{eq4.7} does not hold then we carry on incrementing $\stp_2$.  Split $g_{j}$ as $f^{(2)}_{j,i} + f^{(3)}_{j,i}$, where
\[ f^{(3)}_{j,i}(x) := \sum_{\gamma \in \Omega_{j} \cup S_{j,1} \cup \dots \cup S_{j,i}} \widehat{g_{j}}(\gamma)\overline{\gamma(x)}.\]
If \[ \sup_{\gamma} \sum_{\gamma' \in \gamma + U_{j}} |\widehat{f^{(2)}_{j,i}}(\gamma')| \leq \epsilon^2\] then \textbf{STOP}. Otherwise, there is a $\gamma_{j,i+1}$ with \begin{equation}\label{eq4.81} \sum_{\gamma' \in \gamma_{j,i+1} + U_{j}} |\widehat{f
^{(2)}_{j,i}}(\gamma')| \geq \epsilon^2.\end{equation} Set $\stp = (j,i+1)$, $S_{j,i+1} := (\gamma_{j,i+1} + U_{j}) \setminus (\Omega_{j} \cup S_{j,1} \cup \dots \cup S_{j,i})$. An important point to note is that
\begin{equation}\label{eq4.9} \sum_{\gamma' \in S_{j,i+1}} |\widehat{f^{(2)}_{j,i}}(\gamma')| \geq \epsilon^2:\end{equation} although this is ostensibly a stronger statement than \eqref{eq4.81}, this is merely an illusion since by construction the support
 of $\widehat{f_{j,i}^{(2)}}$ is disjoint from $\Omega_{j} \cup S_{j,1} \cup \dots \cup S_{j,i}$.\vspace{11pt}

\noindent\textit{Claim 2.} The algorithm terminates, and in fact we have the bound $|\Gamma_l| \leq 2\epsilon^{-2}$, independently of $l$.\\
\proof To each element $\gamma_{j,i}$, $j \leq l$ contained in $\Gamma_l$ is associated, by \eqref{eq4.9}, a set $S_{j,i}$ with the property that
\[ \sum_{\gamma \in S_{j,i}} |\widehat{f^{(2)}_{j,i-1}}(\gamma)| \geq \epsilon^2.\] This implies that
\[ 2\sum_{\gamma \in S_{j,i}} |\widehat{f}(\gamma)| \geq \sum_{\gamma \in S_{j,i}} |\widehat{f}(\gamma)||1 - \beta_{j}(\gamma)| = \sum_{\gamma \in S_{j,i}} |\widehat{f^{(2)}_{j,i-1}}(\gamma)| \geq \epsilon^2.\] Recalling that the sets $S_{j,i}$ are disjoint and that $\Vert \widehat{f} \Vert_1 \leq 1$, the claim follows.\vspace{11pt}

\noindent Suppose that the algorithm terminates when $\stp_1 = l$. The whole purpose of having an iteration somewhat more complicated than the one we outlined at the start of the section is that we can beat the crude bound $l \leq 2\epsilon^{-2}$ which follows from Claim 2. \vspace{11pt}

\noindent\textit{Claim 3.} We have $l \leq 2/\epsilon$.\\
\proof Set $S_j = \bigcup_i S_{j,i}$. For any $j$ we have \[ 2\sum_{\gamma \in S_{j}} |\widehat{f}(\gamma)| \geq \sum_{\gamma \in S_{j}} |\widehat{f}(\gamma)||1 - \widehat{\beta_{j}}(\gamma)| = \sum_{\gamma \in S_{j}} |\widehat{g_{j}}(\gamma)| \geq
\epsilon.\]
Since the sets $S_{j}$ are disjoint and $\Vert \widehat{f} \Vert_1 \leq 1$, we do indeed have $l \leq 2/\epsilon$.\vspace{11pt}

\ni To conclude the proof of Lemma \ref{lem4.33}, suppose that $\stp = (l,m)$ at termination. By construction the decomposition
\[ f = f^{(1)}_{l,m} + f^{(2)}_{l,m} + f^{(3)}_{l,m}\] satisfies all of the properties required. It remains to check the bounds claimed. This is a simple matter; using the crude estimate $\eta_{j+1} \geq \epsilon^5\eta'_{j}/2^{11} \geq \epsilon^9\eta_{j}/
2^{21}$, we do indeed have $\eta_l \geq (\epsilon^9/2^{21})^{2/\epsilon}$.\endproof\vspace{11pt}

\ni In order to apply Lemma \ref{lem4.33}, we must be able to exhibit a value of $x$ for which $|f^{(1)}(x)|$ is small. This is where we invoke a kind of ``discrete intermediate value theorem'' of the type used in \S \ref{sec1}. This argument, as the reader will appreciate, is specific to the case $G = \Zp$.

\begin{lemma}\label{lem4.77}
Let $f : \Zp \rightarrow \mathbb{R}$ be a function with $\E_x f(x) = 0$. Let $B = B(\Gamma,\eta)$ be a regular Bohr set with normalised cutoff $\beta$, and suppose that $f^{(1)} = f \ast \beta$. Suppose that $(\epsilon^2 \eta/200d)^d p > 1$. Then
$\min |f^{(1)}(x)| \leq 2\epsilon$.
\end{lemma}
\proof We invoke Lemma \ref{lem4.1}, or rather the proof of it, which tells us that if $\eta' = \epsilon^2 \eta/200 d$ and if $t \in B(\Gamma,\eta')$ then
\begin{equation}\label{eq4.21} |f^{(1)}(x + t) - f^{(1)}(x)| \leq 2\epsilon\end{equation}
for any $x$. The conditions of the lemma, together with the simple bound of Lemma \ref{lem4.00} (i), imply that $|B(\Gamma,\eta')| > 1$, so there is some $t \in B(\Gamma,\eta')$, $t \neq 0$. The elements $\{0,t,2t,3t,\dots, (p-1)t\}$ are in fact just the
elements of $\Zp$, listed once each. Since $\E_x f^{(1)}(x) = 0$, there must be a value of $j$ for which
$f^{(1)}(jt)$ and $f^{(1)}((j+1)t)$ have opposing signs. For this value of $j$,
\eqref{eq4.21} guarantees that $|f^{(1)}(jt)| \leq 2\epsilon$.\endproof
\vspace{11pt}

\ni\textit{Proof of Theorem \ref{mainthm2}.} Lemma \ref{lem4.33} tells us that there is a set $\Gamma$, $d := |\Gamma| \leq 2/\epsilon^2$, and an $\eta \geq (\epsilon^9/2^{21})^{2/\epsilon}$ for which the hypotheses of Lemma \ref{lem4.14} are satisfied.
If in addition $\min_x |f^{(1)}(x)| \leq 2\epsilon$, then the conclusion of that lemma
implies that $\min_x |f(x)| \leq 8\epsilon$. By Lemma \ref{lem4.77}, this will be the case if $(\epsilon^2 \eta/200d)^d p > 1$, which will be satisfied if $\epsilon
=C(\log p/\log \log p)^{1/3}$ for suitable $C$.\endproof\vs

\ni It is a simple matter to deduce Theorem \ref{mainthm1} from Theorem \ref{mainthm2}. Recall that $S(p)$ is the smallest possible value of $\Vert \widehat{1_A} \Vert_1$, where $A \subseteq \Zp$ is a set with cardinality $\lfloor p/2 \rfloor$. \vs

\ni\textit{Proof of Theorem \ref{mainthm1}.} Set $f(x) := (21_A(x) - 1)/\Vert (21_A - 1)^{\wedge} \Vert_1$, so that $\Vert \widehat{f} \Vert_1 = 1$, and apply Theorem \ref{mainthm2}. One obtains
\[ \Vert (21_A - 1)^{\wedge} \Vert_1^{-1} = \min_x |f(x)| = O(\log \log p/\log p)^{1/3}.\]
Now
\[ (21_A - 1)^{\wedge}(r) = \left\{ \begin{array}{ll} 2\widehat{1_A}(r) & \mbox{if $r \neq 0$} \\ O(1/p) & \mbox{if $r = 0$},   \end{array}  \right.\]
which means that
\[ \Vert (21_A - 1)^{\wedge} \Vert_1 \geq 2 \Vert \widehat{1_A} \Vert_1 + O(1).\] The result follows immediately.\endproof\vs

\ni To conclude this section we give the proof of Theorem \ref{mainthm3.1} and hence, by the remarks in \S \ref{sec1}, answer the question of Gowers. \vs

\ni\textit{Proof of Theorem \ref{mainthm3.1}.} Consider the functions
\[ f_0(x) := f(x) - \E f,\]
\[ g_0(x) := g(x) - \E g\]
and
\[ h(x) := f \ast g - \E f \E g.\]
It is easy to see that $\E f_0 = \E g_0 =0$, that $\Vert f_0 \Vert_2, \Vert g_0 \Vert_2 \leq 1$, and that
\[ h = f_0 \ast g_0.\]
In particular we have, by the Cauchy-Schwarz inequality and Parseval's identity, that
\[ \Vert \widehat{h} \Vert_1 \leq \Vert f_0 \Vert_2\Vert g_0 \Vert_2 \leq 1\]
(this is a special case of Young's inequality).
The result follows immediately from Theorem \ref{mainthm2}.\endproof

\section{The proof of Theorem \ref{mainthm3.2}}\label{sec6}
\begin{lemma}\label{lem5.10} For any prime $p$ there is a function
$F : \Zp \rightarrow \mathbb{R}$ such that $\E F=0$,
$\Vert \widehat{F} \Vert_1\leq 1$, $\widehat{F}(r) = \widehat{F}(-r) \geq 0$
for all $r \in \Zp$, and $\min_x |F(x)|\gg(\log p)^{-1}$.
\end{lemma}
\proof We may, of course, assume that $p$ is sufficiently large. Consider to begin with the very simple function
\[ g_1(x) := \left\{ \begin{array}{ll} 1 & \mbox{if $|x/p| < 1/4$} \\ -1 & \mbox{otherwise}.  \end{array}\right.\]
We have $\E g_1 = O(1/p)$, whilst the Fourier transform $\widehat{g}_1(r)$ satisfies $\widehat{g}_1(r) = \widehat{g}_1(-r)$ and
\begin{equation}\label{eq5.1a} |\widehat{g}_1(r)| \leq \frac{1}{p|\sin(\pi r/p)|} \leq \frac{1}{2|r|} \end{equation}
for $|r| < p/2$. For suitable $C$, the function $F_1(x) := g_1(x)/C\log p$ satisfies all the requirements of the lemma except for the condition that $\widehat{F}_1(r)$ be non-negative (and the condition that $\E F_1 = 0$, though this can be achieved by a
very small perturbation of $F_1$).\vs

\ni To construct a function with non-negative Fourier transform, we begin by considering an auxilliary function on $\mathbb{R}/\Z$. Define
\[ f_1(\theta) := -2 \log |2 \sin \pi \theta|\] for $\theta \neq 0$, and $f_1(0) = 0$.
Note that $f \in L^p$ for all $p \in [1,\infty)$. We claim that the Fourier transform of $f_1$ is given by
\begin{equation}\label{eq5.1b}
\widehat{f}_1(m) = \left\{ \begin{array}{ll} 1/|m| & \mbox{if $m \in \Z \setminus \{0\}$} \\ 0 & \mbox{otherwise}.  \end{array}\right.\end{equation}
Heuristically, this can be seen by observing the Taylor expansion
\[ -\log(1 - z) = z + \frac{z^2}{2} + \frac{z^3}{3} + \dots,\]
setting $z = e^{2\pi i \theta}$ and $z = e^{-2\pi i \theta}$, and adding. To justify the assertion rigourously, one could proceed as follows. Using Cauchy's formula for derivatives, one confirms that
\begin{equation}\label{eq5.1c} \frac{1}{2\pi i} \int_{\Gamma(\epsilon)} \frac{\log (1 - z)}{z^{n+1}} \, dz = -\frac{1}{n}\end{equation} for $n \geq 1$, where $\Gamma(\epsilon)$ is the contour consisting of the circle $|z| = 1$ indented to the left near $z
 = 1$ using a semicircle of radius $\epsilon$. The contribution from the semicircle of radius $\epsilon$ is $O(\epsilon \log(1/\epsilon))$, which tends to 0 as $\epsilon \rightarrow 0$. The remainder of the integral is
\[ \int_{\theta \in \mathbb{R}/\Z : |\theta| \geq \epsilon} \log(1 - e^{2\pi i \theta}) e^{-2\pi i n \theta} \, d\theta.\] This differs from the integral over all $\theta \in \mathbb{R}/\Z$ by an error of
\[ O \big( \int^{\epsilon}_0 |\log \theta | \, d\theta\big),\] which is $O(\epsilon \log(1/\epsilon))$. Letting $\epsilon \rightarrow 0$ in \eqref{eq5.1c}, then, we see that
\begin{equation}\label{eq5.1d} \int_{\theta \in \mathbb{R}/\Z} \log(1 - e^{2\pi i \theta}) e^{-2\pi i n\theta}\, d\theta = - \frac{1}{n}\end{equation} for $n \geq 1$. By a similar argument, this integral vanishes when $n \leq 0$. Furthermore by an almost
identical argument one may confirm that
\[ \int_{\mathbb{R}/\Z} \log(1 - e^{-2\pi i \theta}) e^{-2\pi i n\theta} = \left\{ \begin{array}{ll} 1/n & \mbox{if $n \leq -1$} \\ 0 & \mbox{otherwise}.\end{array}\right.\]
Adding this to \eqref{eq5.1d} establishes the claim \eqref{eq5.1b}.\vs

\ni Now $f_1$ is useless to us as it stands, since $\widehat{f}_1$ does not lie in $l^1(\mathbb{Z})$. We must also transfer $f_1$ to a function defined on $\Z/p\Z$. We modify $f_1$ by defining
\begin{equation}\label{eq5.1h} f_2 := f_1 \ast \chi \ast \chi,\end{equation}
where
\[ \chi(\theta) = \frac{Cp}{2} 1_{|\theta| \leq 1/Cp}(\theta)\]
for some large constant $C$ to be specified later. Note that
\begin{equation}\label{eq5.1e} \widehat{f}_2(m) = \widehat{f}_1(m) \big(\frac{Cp\sin 2\pi m/Cp}{2\pi m}\big)^2,\end{equation} where the bracketed expression is to be interpreted as 1 if $m = 0$.
Consider the function $g_2 : \Z/p\Z \rightarrow \mathbb{R}$ defined by
\begin{equation}\label{eq5.1o} g_2(x) = \sum_{m \in \Z} \widehat{f}_2(m) e^{2\pi i m x/p}.\end{equation}
Since $\widehat{f}_2 \in l^1(\mathbb{Z})$, the series converges uniformly and
\[ g_2(x) = f_2(x/p)\] provided that $x \neq 0$. \vs

\ni Now we have
\[ |f'_1(\theta)| = |2\pi \cot \pi \theta| \leq \frac{2}{|\theta|}\] for all $0 < |\theta| \leq 1/2$.
From the definition \eqref{eq5.1h}, it follows that
\begin{equation}\label{eq5.1i} |g_2(x) - f_1(x/p)| \leq \frac{2}{Cp} \sup_{|\theta - x/p| \leq 2/pC} |f'_1(\theta)| \leq 8/C\end{equation} provided that $x \neq 0$ and $C \geq 4$.\vs

\ni All we need to know about $g_2(0)$ is that it is large, to which end the bound
\begin{equation}\label{eq5.1j} g_2(0) \geq g_2(1)\end{equation} will be quite sufficient.\vs

\ni The Fourier transform $\widehat{g}_2(r)$ on $\Z/p\Z$ is given by
\begin{equation}\label{eq5.1f} \widehat{g}_2(r) =  \sum_{\substack{m \in \Z \\ m \equiv r \mdsub{p}}} \widehat{f}_2(m) \end{equation} (this is an instance of the Poisson summation formula, but in this case it follows from the inversion formula
\[ g_2(x) = \sum_r \widehat{g}_2(r) e^{2\pi i rx/p}\]
upon comparing with \eqref{eq5.1o} and recalling the uniqueness of Fourier expansion). In particular $\widehat{g}_2(r) = \widehat{g}_2(-r)$. Now $\widehat{f}_2(m)$ is real and non-negative, and hence so is $\widehat{g}_2(r)$ and we have
\[ \widehat{g}_2(r) \geq \widehat{f}_2(\overline{r}),\]
where $\overline{r}$ is the unique integer with $|\overline{r}| < p/2$ and $\overline{r} \equiv r \md{p}$. Using the inequality $\sin \theta/\theta \geq 1 - \theta^2/6$ and \eqref{eq5.1b}, we thus have
\begin{equation}\label{eq5.1ff} \widehat{g}_2(r) \geq \widehat{f}_1(\overline{r}) \big( 1 - \frac{2\pi^2 |\overline{r}|^2}{3C^2 p^2}\big)^2 \geq \frac{1}{|\overline{r}|}\big(1 - \frac{4}{C^2}\big)\end{equation} for $r \neq 0$.\vs

\ni In the other direction, we have the estimate
\begin{equation}\label{eq5.1g} \Vert \widehat{g}_2 \Vert_1 = \sum_m \widehat{f}_2(m) = \sum_{\substack{m \in \Z \\ m \neq 0}}  \frac{1}{|m|}  \big(  \frac{Cp  \sin  2\pi  m/Cp}{2\pi  m}\big)^2  \leq  2\sum_{m=1}^{\infty}  \min(\frac{1}{m},  \frac{C^2
p^2}{4\pi^2m^3}) \ll \log p. \end{equation}
\ni Now define
\[ g_3 := g_1 + g_2.\]
From \eqref{eq5.1a} and \eqref{eq5.1ff}, we see that $\widehat{g}_3(r) = \widehat{g}_3(-r)$, and that $\widehat{g}_3(r) \geq 0$ for all $r \neq 0$, provided that $C^2 \geq 8$. From \eqref{eq5.1a} and \eqref{eq5.1g} we see that
\begin{equation}\label{eq5.1m} \Vert \widehat{g}_3 \Vert_1 \ll \log p.\end{equation} We claim that if $C \geq 1000$ then \[ |g_3(x)| \geq 1/4 \qquad \mbox{for all $x$}.\]
To see this, observe from \eqref{eq5.1i} that if $1/4 < |x/p| \leq 1/2$ then
\[ g_3(x) = -1 + g_2(x) \leq -1 + f_1(x/p) + 8/C \leq -1 - \log 2 + 8\pi/C \leq -1.\]
If $0 < |x/p| < 1/4$ then we have
\[ g_3(x) = 1 + g_2(x) \geq 1 + f_1(x/p) - 8/C \geq 1 - \log 2 - 8\pi/C \geq 1/4.\]
Finally from \eqref{eq5.1j} we have
\[ g_3(0) = 1 + g_2(0) \geq 1 + g_2(1/p) \geq 1 + f_1(1/p) - 8/C \geq 1/4.\]
This proves the claim.\vs

\ni The function $F$ (whose construction is our goal) will essentially be a rescaled version of $g_3$, but it is required to satisfy $\E F = 0$. To this end we note that
\[ \E g_1 = O(1/p),\]
whilst, from \eqref{eq5.1f}, we have
\[ \E g_2 = \sum_{m \equiv 0 \mdsub{p}} \widehat{f}_2(m) \leq 2\sum_{\substack{m =1 \\ m\equiv 0 \mdsub{p}}}^{\infty} \min(\frac{1}{m}, \frac{C^2p^2}{4\pi^2 m^3}) = O(1/p).\]
Thus $\E g_3 = O(1/p)$, and so if we define
\[ g_4 := g_3 - \E g_3\]
then we still have
\[ |g_4(x)| \geq 1/4 + O(1/p) \geq 1/5\]
for $p$ sufficiently large. Of course, since $\widehat{g}_4(r) = \widehat{g}_3(r)$ for $r \neq 0$, we also have $\widehat{g}_4(r) = \widehat{g}_4(-r) \geq 0$.\vs

\ni At last we may define $F$ by
\[ F(x) := g_4(x)/\Vert \widehat{g}_4 \Vert_1.\] The asserted properties of $F$ are immediate from the definition, the facts we have assembled about $g_4$, and \eqref{eq5.1m}.\endproof\vs

\ni We conclude this section by proving Theorem \ref{mainthm3.2}, which asserted the existence of a function $f : \Zp \rightarrow \mathbb{R}$ with $\E f = 0$, $\Vert f \Vert_2 \leq 1$, and $|f \ast f^{\circ}(x)| \geq c/\log p$ for all $x$.\vs

\ni\textit{Proof of Theorem \ref{mainthm3.2}.} Let $F$ be the function
constructed in Lemma \ref{lem5.10}. Take the function $f$
with $\widehat{f}=\sqrt{\widehat{F}}$. Then
\[\E f = \widehat{f}(0) = \sqrt{\widehat{F}(0)}
= \sqrt{\E F} = 0,\]
\[\E f^2 = \sum_r \widehat{f}(r)^2 = \Vert \widehat{F} \Vert_1\leq 1,\]
and
\[ f \ast f^{\circ} = F.\]
The theorem follows. \endproof

\section{An example for Gowers' question}\label{sec7}

\ni In this section we prove Theorem \ref{mainthm3.3}, which asserted that there is a set $A \subseteq \Zp$ with $|A| = \lfloor p/2\rfloor$ and
\[ |A \cap (A + x) - p/4| \geq c/\log p\log\log p\] for all $x$. Noting that
\[ |A \cap (A + x) - p/4| = p f \ast f^{\circ}(x),\]
where $f := 1_A - 1/2$, we see that this is a matter of finding a function $f$ in Theorem \ref{mainthm3.2} which somehow ``resembles a set''. The function $f$ constructed at the end of \S \ref{sec6} need not have this property. Note, however, that we could have exercised considerable freedom in our choice of $f$: any function such that $\widehat{f}(r) = \xi_r \sqrt{\widehat{F}(r)}$, where $\xi_r \in \{-1,1\}$, would have done. By choosing the signs $\xi_r$ at random we may force $f$ to behave well in $L^{\infty}$, and hence to resemble a set. The rest of the section is devoted to the details of such an argument.\vs

\ni Suppose, then, that $\xi = \{\xi_r\in\{\pm1\},\, 0\leq r < p/2\}$ is a sequence of independent random signs and that $f=f_{\xi}$ is the function whose Fourier transform satisfies $\widehat{f}(r) = \xi_r \sqrt{\widehat{F}(r)}$, $0 \leq r < p/2$, and $\widehat{f}(r) = \widehat{f}(-r)$.\vs

\ni Let us reiterate the properties enjoyed by $f_{\xi}$: It satisfies $\E f_{\xi} = 0$, $\Vert f_{\xi} \Vert_2 \leq 1$, and $|f_{\xi} \ast f_{\xi}^{\circ}(x)| \gg 1/\log p$ for all $x$.

\begin{lemma}\label{lem6.1} There is an absolute constant $c$ and a choice of signs $\xi$ such that
\begin{equation}\label{largedev}
\E_{x \in \Zp} \exp\left(cf_{\xi}(x)^2\right) \ll 1.
\end{equation}
\end{lemma}

\ni\textit{Proof.} Let $(a_r)_{0 \leq r < p/2}$ be any sequence of real numbers. It is well-known (a simple instance of Khintchine's inequality, see for example \cite{kah}) that
\[ \E_{\xi}| \sum_r a_r \xi_r |^{2k} \leq \big(C_0 k \sum_r a_r^2\big)^k \]
for some absolute $C_0$. It follows that for any $c \geq 0$ we have
\begin{equation}\label{eq3.394} \E_{\xi} \exp \big( c \big(\sum_r a_r \xi_r \big)^2\big) \leq \sum_{k \geq 0}\frac{1}{k!} \big(cC_0 k \sum_r a_r^2 \big)^k.\end{equation}
Now for any $x \in \Zp$ we have
\[ f_\xi(x)=\sum_{0 < r < p/2}2\xi_r\sqrt{\widehat{F}(r)}\cos(2\pi xr/p).\]
For each fixed $x$ this has the form
\[ \sum_r a_r \xi_r\]
where
\[ a_r = 2\sqrt{\widehat{F}(r)} \cos (2\pi x r/p).\] Since $\Vert \widehat{F} \Vert_1 \leq 1$, we have
\[ \sum_r a_r^2 \leq 4.\]
In view of \eqref{eq3.394} this means that
\[ \E_{\xi} \exp(cf_{\xi}(x)^2) \leq \sum_{k \geq 0} \frac{1}{k!}(4cC_0k)^k.\]
By Stirling's formula this series converges provided that $4c C_0 < 1/e$, an inequality which can certainly be achieved by selecting $c$ sufficiently small.\vs

\ni Averaging over $x$, we see that
\[ \E_x \E_{\xi} \exp(c f_{\xi}(x)^2) \ll 1.\] In particular there is at least one choice of the signs $\xi$ for which the claimed inequality is true.\endproof\vs

\ni The next lemma is a technical fact used in the proof of Lemma \ref{lem6.3} below.

\begin{lemma}\label{lem6.2} Let $f_1, f_2 : \Zp \rightarrow \mathbb{R}$ be two functions. Then the following inequality holds:
\begin{equation}\label{stab}
\Vert f_1 \ast f_1^{\circ}- f_2 \ast f_2^{\circ}\Vert_{\infty} \leq \|f_1-f_2\|_2
\left(\|f_1\|_2+\|f_2\|_2\right).
\end{equation}
\end{lemma}

\ni\textit{Proof.} By the Cauchy--Schwarz inequality we have
\[ |f*g(x)| \leq \|f\|_2\|g\|_2\] for any two functions $f,g$ and for any $x \in \Z/p\Z$.
Setting $g:=f_1-f_2$ we have
$$|(f_1 \ast f_1^{\circ})(x)-(f_2 \ast f_2^{\circ})(x)| \leq |f_2*g^{\circ}(x)| +
|g* f_1^{\circ}(x)| \leq \|f_2\|_2\|g\|_2 + \|f_1\|_2\|g\|_2,$$
as required. \endproof\vs

\ni The next lemma is a stronger version of Theorem \ref{mainthm3.2}.

\begin{lemma}\label{lem6.3} There is a function $f : \Z/p\Z \rightarrow \R$ such that $\E f = 0$, $\Vert f \Vert_2 \ll 1$, $|f \ast f^{\circ}(x)| \gg 1/\log p$ for all $x$, and with the additional property that $\Vert f \Vert_{\infty} \ll \sqrt{\log \log
p}$.
\end{lemma}

\ni\textit{Proof.} Take a function $f_{\xi}$ satisfying the conclusion of Lemma \ref{lem6.1}, that is to say for which the inequality
\begin{equation}\label{eq4.72} \E_{x} \exp(c f_{\xi}(x)^2) \ll 1\end{equation} is satisfied. Let $C$ be a large absolute constant, and consider the set
\[ A := \{x : |f_{\xi}(x)| \geq C\sqrt{\log \log p}\}.\]
It follows from \eqref{eq4.72} that, for a suitable choice of $C$, we have
\begin{equation}\label{exceptional-size} \E 1_A \ll (\log p)^{-10}.\end{equation}
Choose an arbitrary set $A'$ such that $A \subseteq A'$ and
\begin{equation}\label{exact-size} (\log p)^{-10} \ll \E 1_{A'} \ll (\log p)^{-10}.\end{equation}
Define $f$ by
\[ f(x) := \left\{ \begin{array}{ll} f_{\xi}(x) & \mbox{if $x \notin A'$}; \\ \kappa & \mbox{if $x \in A'$},\end{array}\right.\]
where $\kappa$ is selected so that $\E f = 0$. \vs

\ni Observe that
\[ |\kappa| = |\E_{x \in A'} f_{\xi}(x)| \leq \E_{x \in A'} |f_{\xi}(x)|.\]
From \eqref{eq4.72} one may deduce without difficulty that
\[ \E_{x \in S} |f_{\xi}(x)| \ll \big( \log (1/\E1_S) \big)^{1/2}\]
for any nonempty set $S \subseteq \Z/p\Z$. Hence, in view of \eqref{exact-size}, we have the bound
\[ |\kappa| \ll \sqrt{\log \log p}.\]
This immediately implies that $\Vert f \Vert_{\infty} \ll \sqrt{\log \log p}$. \vs

\ni To obtain the upper bound on $\Vert f \Vert_2$, it suffices to note that
\[ \Vert f - f_{\xi} \Vert_2 = \Vert (f_{\xi} - \kappa) 1_A \Vert_2 \ll \sqrt{\log p} (\E 1_A)^{1/2} \ll (\log p)^{-4}.\]
Furthermore, this inequality together with Lemma \ref{lem6.2} implies that
\[ \Vert f \ast f^{\circ} - f_{\xi} \ast f_{\xi}^{\circ}\Vert_{\infty} \ll (\log p)^{-4}.\]
In view of the fact that 
\[ |f_{\xi} \ast f_{\xi}^{\circ}(x)| \gg 1/\log p\] for all $x  \in \Z/p\Z$,
we obtain the desired lower bound
\begin{equation}\boxeq |f \ast f^{\circ}(x)| \gg 1/\log p.\end{equation}

\ni\textit{Remark.} A well-known conjecture of Koml\'os implies that if $(a_{ij})_{1 \leq i,j \leq n}$ are any scalars then there is a choice of signs $\xi_j$ so that
\[ \big| \sum_{j} a_{ij} \xi_j \big| \leq C \big( \sum_{j} a_{ij}^2 \big)^{1/2},\]
for some absolute constant $C$. This conjecture, if true, would allow us to replace the $\sqrt{\log \log p}$ appearing in Lemma \ref{lem6.3} by an absolute constant, and to remove all subsequent appearances of $\log \log p$. The best result towards Koml\'os' conjecture is due to Banaszczyk \cite{bana}; it may be that his ideas may be applied in our context to reduce the $\log \log p$ factors to something a little smaller, but it does not seem worth the effort of pursuing this line of inquiry.\vs

\ni Using a simple linear transformation, Lemma \ref{lem6.3} has the following consequence.

\begin{lemma}\label{lem6.4} There exists a function $f: \Zp \rightarrow [0,1]$
such that $\E f=1/2$ and 
\begin{equation}\label{lem6.4bd} |f \ast f^{\circ}(x) - 1/4| \gg 1/\log p\log\log p \end{equation} for all $x \in \Z/p\Z$.\endproof
\end{lemma}

\ni To complete the proof Theorem \ref{mainthm3.3} we use probabilistic arguments again, in a fairly standard manner. We choose a set $B$ by selecting each $x \in \Z/p\Z$ to lie in $B$ independently at random with probability $f(x)$, and then show that with high probability $|B|$ is close to $p/2$ and $1_B \ast 1_B^{\circ}(x)$ is close to $f \ast f^{\circ}(x)$ whenever $x \neq 0$. We then modify $B$ by adding or deleting a small number of elements to give a set $A$ with $|A| = \lfloor p/2 \rfloor$.\vs

\ni For each $x \in \Z/p\Z$ we write $X_x$ for the random variable which is $1$ with probability $f(x)$ and $0$ with probability $1 - f(x)$. To prove the statements we need, it is natural to use one of the standard large-deviation type inequalities. Amongst inequalities of this sort, one particularly suited to our purpose is Hoeffding's inequality (see, for example, \cite{green-deviation-notes}).

\begin{lemma}[Hoeffding's inequality]
Suppose that $Y_1,\dots, Y_m$ are independent random variables such that $|Y_i| \leq 1$ for each $i$. Write $Y := m^{-1}(Y_1 + \dots + Y_m)$ and $\mu := \E Y$. Then for any $t \geq 0$ we have
\[ \mathbb{P}(|Y - \mu| \geq t) \leq 2 e^{-mt^2/2}.\]
\end{lemma}
\ni It follows immediately from this, applied with $m = p$ and $Y_x = X_x$, that 
\begin{equation}\label{b-size} \mathbb{P}(| |B| - p/2 | \geq p^{2/3}) \ll e^{-p^{1/3}/2}.\end{equation}
Fix $x \in \Z/p\Z$, $x \neq 0$, and look at the expression
\[ 1_B \ast 1_B^{\circ}(x) = p^{-1}\sum_y X_y X_{y - x}.\]
Writing 
\[ Z_y := X_y X_{y - x},\]
we have 
\[ p^{-1}\sum_y \E Z_y = f \ast f^{\circ}(x).\]
Hoeffding's inequality does not immediately apply to the family $(Z_y)_{y \in \Z/p\Z}$, since these random variables are not independent. However we may partition the family into three subfamilies $(Z_y)_{y \in S_i}$, $i = 1,2,3$, where $|S_i| \geq p/5$, and within each family the random variables $Z_y$ are jointly independent. This may be achieved by forming the graph on vertex set $\Z/p\Z$ in which verex $i$ is joined to vertex $j$ if $i - j = x$; this is a cycle of length $p$, and we partition its vertex set into 3 classes each having size at least $p/5$.\vs

\ni Hoeffding's inequality now applies, and we have
\[ \mathbb{P}\big( \big| |S_i|^{-1}\sum_{y \in S_i} Z_y - |S_i|^{-1}\sum_{y \in S_i} \E Z_y \big| \geq t \big) \leq 2 e^{-p t^2/10} \]for $i = 1,2,3$. It follows using the triangle inequality that
\[ \mathbb{P}\big( \big| p^{-1} \sum_{y \in \Z/p\Z} Z_y - p^{-1} \sum_{y \in \Z/p\Z} \E Z_y \big| \geq t \big) \leq 6 e^{-p t^2/10},\] that is to say
\[ \mathbb{P}\big( \big| (1_B \ast 1_B^{\circ})(x) - (f \ast f^{\circ})(x) \big| \geq t \big) \leq 6e^{-p t^2/10}.\]
Taking $t = p^{-1/3}$, we see from this and \eqref{b-size} that with probability close to 1 we have both
\begin{equation}\label{eq-19a} \big| |B| - p/2 \big| \leq p^{2/3}\end{equation}
and
\begin{equation}\label{eq-19b} \big| 1_B \ast 1_B^{\circ}(x) - f \ast f^{\circ}(x) \big| \leq p^{-1/3}\end{equation} for all $x \neq 0$. Fix a specific set $B$ verifying \eqref{eq-19a} and \eqref{eq-19b}.
Form a set $A$ by adding to or deleting from $B$ arbitrary elements so that $|A| = \lfloor p/2\rfloor$. In view of \eqref{eq-19a} we have
\[ \Vert 1_A - 1_B \Vert_2 \leq 2p^{-1/6},\] and hence by Lemma \ref{lem6.2} we infer that
\[ \Vert 1_A \ast 1_A^{\circ} - 1_B \ast 1_B^{\circ} \Vert_{\infty} \leq 4p^{-1/6}.\] It follows from this, \eqref{eq-19b} and \eqref{lem6.4bd} that
\[ | 1_A \ast 1_A^{\circ}(x) - 1/4 | \gg 1/\log p \log\log p\] for all $x \neq 0$. This inequality is manifestly true for $x = 0$ as well, and this completes the proof of Theorem \ref{mainthm3.3}.\endproof

\section{A recent development}

\ni We have recently learnt that Sanders \cite{sanders} has combined some of our methods with some new ideas, and has been able to improve the bound in Theorem \ref{mainthm2} to $(\log p)^{-1/2 + \epsilon}$.

\end{document}